\documentclass[10pt]{amsart}

\usepackage{amsmath,amsfonts,amscd,amssymb,epsf}
\newtheorem{theorem}{Theorem}[section]
\newtheorem{proposition}[theorem]{Proposition}
\newtheorem{lemma}[theorem]{Lemma}

\theoremstyle{definition}

\theoremstyle{remark} 
\numberwithin{equation}{section}
\DeclareMathOperator{\PSU}{PSU}

\newcommand{\C}{\mathbb{C}}
\newcommand{\Del}{\mathbb{D}}

\newcommand{\R}{\mathbb{R}}
\newcommand{\bk}{\backslash}
\newcommand{\pa}{\partial}

\newcommand{\ov}{\overline}

\begin{document}

\title{The Weil--Petersson geometry of the moduli space of Riemann surfaces}
\author{Lee-Peng Teo}\email{lpteo@mmu.edu.my}

\address{Faculty of Information Technology, Multimedia University,
Jalan Multimedia, Cyberjaya, 63100, Selangor Darul Ehsan, Malaysia.}

 \maketitle

\begin{abstract}In \cite{Huang1}, Z. Huang showed that in the thick part of the moduli
space $\mathcal{M}_g$ of compact Riemann surfaces of genus $g$, the
sectional curvature of the Weil--Petersson metric is bounded below
by a constant depending on injectivity radius, but independent of
the genus $g$. In this article, we prove this result by a different
method. We also show that the same result holds for Ricci curvature.
For the universal Teichm\"uller space equipped with Hilbert
structure induced by Weil--Petersson metric, we prove that its
sectional curvature is bounded below by a universal constant.

\end{abstract}

\section{Introduction}There have been a lots of studies on the geometry of Weil--Petersson (WP) metric on
moduli spaces of Riemann surfaces, especially regarding its
curvature properties \cite{Ahlfors, Royden, Wolpert1, Tromba,
Schumacher, Jost, McMullen, Wolpert2, Yamada, Wilson, Liu_Sun_Yau_1,
Liu_Sun_Yau_2, Huang2, Huang1}. In a pioneering work, Ahlfors
\cite{Ahlfors} showed that the Ricci, holomorphic sectional and
scalar curvatures of the WP metric are all negative. Later, Royden
\cite{Royden} showed that the holomorphic sectional curvature is
bounded above by a negative constant, and he conjectured that on the
moduli space $\mathcal{M}_g$ of compact Riemann surfaces of genus
$g$, this constant is equal to $\frac{-1}{2\pi(g-1)}$. By deriving
more compact expressions for the Riemann tensors of the WP metric,
Wolpert \cite{Wolpert1} verified Royden's conjecture. He also showed
that the Ricci curvature is  bounded above by
$\frac{-1}{2\pi(g-1)}$, and the scalar curvature is bounded above by
$-\frac{3(3g-2)}{4\pi}$. In the communications between Wolpert and
Tromba and between Wolpert and Royden, it was
  proved that the sectional curvature of the WP metric is also
negative. Detailed proof of this result was given by Wolpert in
\cite{Wolpert1} and by Tromba in \cite{Tromba}. Regarding the upper
bound, it was proved in \cite{Huang2} that the sectional curvature
does not have a negative upper bound.

Lower bounds of the curvatures are less considered. There are some
results obtained by \cite{Schumacher, Trapani, Huang3, Huang1}. The
results of \cite{Schumacher, Trapani, Huang3} showed that the
sectional curvature  is not bounded below on the moduli space
$\mathcal{M}_g$. Therefore, attention must  be shifted to find lower
bounds of the sectional curvature on compact subsets of the moduli
space $\mathcal{M}_g$. This problem was  studied by Huang in
\cite{Huang1}. To describe his result in more detail, we need to
introduce some notions first. A point on the moduli space
$\mathcal{M}_g$ can be considered as a compact Riemann surface $X$
of genus $g$ endowed with a unique metric of constant curvature
$-1$, called the hyperbolic metric. The injectivity radius of $X$ at
a point $z\in X$, $inj(X;z)$, is defined as the supremum of all $r$
for each the open set $U_z^r=\left\{ w\in X\,:\, d(z,w)<r\right\}$
is isometric to a disc. The injectivity radius of $X$, $inj(X)$, is
defined to be the infimum of $inj(X;z),\, z\in X$. By a well--known
result, $inj(X)$ is equal to one half of the length of the shortest
closed geodesic of $X$. Given a positive constant $r_0$, the thick
part of moduli space $\mathcal{M}_g$ (with respect to $r_0$), is
defined as the subset of $\mathcal{M}_g$ consisting of those points
where the injectivity radius of the corresponding Riemann surfaces
is greater than $r_0$. In \cite{Huang1}, Huang showed that on the
thick part of the moduli space $\mathcal{M}_g$, the holomorphic
sectional  and sectional curvatures of the WP metric are both
bounded below by negative constants $-C_1$ and $-C_2$ depending on
$r_0$, but independent of genus $g$. As a result, the Ricci
  and scalar curvatures are bounded below by $-C_3 g$ and
$-C_4g^2$ respectively, where $C_3$ and $C_4$ are two positive
constants depending on $r_0$, but independent of $g$.
 The main
tool used by Huang is the analysis of harmonic maps between
hyperbolic surfaces. In the present article, we are going to give a
different proof of Huang's result without using harmonic maps.
Moreover, we are going to improve the bounds $-C_3g$ and $-C_4g^2$
for Ricci and scalar curvatures to $-C_3$ and $-C_4g$ respectively.
Explicit dependence of the constants $C_1, C_2, C_3, C_4$ on the
injectivity radius $r_0$ is given.

In \cite{TT}, we have defined a Hilbert structure on the universal
Teichm\"uller space $T(1)$, so that the Weil--Petersson metric is a
well-defined metric on $T(1)$. We have also obtained an explicit
formula for the Riemann curvature tensor of the WP metric, which is
generalization of the result of Wolpert \cite{Wolpert1}. In this
article, we are going to show that  the sectional curvature of the
WP metric on $T(1)$ is bounded below by a universal negative
constant. We also show that it does not have a negative upper bound.

The layout of this article is as follows. In Section 2, we review
some necessary facts. In Section 3, we obtain the lower bounds of
the curvatures of the WP metric on the moduli space $\mathcal{M}_g$
as a function of injectivity radius. In Section 4, we find the lower
bound of the sectional curvature of the WP metric on the universal
Teichm\"uller space.

\section{Background}
In this section, we present some necessary facts. Let $T(X)$ and
$\mathcal{M}(X)$ be respectively the Teichm\"uller space and moduli
space of a compact  Riemann surface $X$ of genus $g$, where $g\geq
2$. The Teichm\"uller space $T(X)$ has a complex analytic model
described as follows. Let $\Del$ and $\Del^*$ be respectively  the
unit disc and its exterior. There is a Fuchsian group $\Gamma\in
\PSU(1,1)$ such that the quotient of the unit disc $\Del$ by the
action of $\Gamma$ is $X$, i.e., $X\simeq \Gamma\bk\Del$. The space
of bounded Beltrami differentials on $X$ can be identified with the
space of bounded $\Gamma$-automorphic $(-1,1)$ differentials on
$\Del$, denoted by $\mathcal{A}^{-1,1}(\Del, \Gamma)$, which
consists of bounded functions $\mu$ on $\Del$ satisfying
\begin{align*}
\mu(\gamma (z)) \frac{\ov{\gamma'(z)}}{\gamma'(z)}=\mu(z).
\end{align*}Let $\mathcal{B}^{-1,1}(\Del,\Gamma)$ be the unit ball
of $\mathcal{A}^{-1,1}(\Del,\Gamma)$ with respect to the sup--norm:
\begin{align*}
\Vert\mu\Vert_{\infty}=\sup_{z\in\Del} |\mu(z)|.
\end{align*}Given a Beltrami differential $\mu\in
\mathcal{B}^{-1,1}(\Del,\Gamma)$, extend it to
$\hat{\C}:=\C\cup\{\infty\}$ by reflection:
\begin{align}\label{eq1}
\mu(z)=\ov{\mu\left(\frac{1}{\bar{z}}\right)}\frac{z^2}{\bar{z}^2},\hspace{1cm}z\in
\Del^*.
\end{align}There is a unique quasiconformal mapping
$w_{\mu}:\hat{\C}\rightarrow\hat{\C}$ which fixes the points $-1,
-i, 1$ and satisfies the Beltrami equation $(w_{\mu})_{\bar{z}}=\mu
(w_{\mu})_z$. The conjugation of $\Gamma$ by $w_{\mu}$,
$\Gamma_{\mu}=w_{\mu}\circ \Gamma\circ w_{\mu}^{-1}$, is again a
Fuchsian group. The corresponding quotient surface
$X_{\mu}=\Gamma_{\mu}\bk\Del$ is a Riemann surface having the same
type as $X$, but with different complex structure. Define an
equivalence relation
 on $\mathcal{B}^{-1,1}(\Del,\Gamma)$ so that $\mu\sim\nu$ if  and
 only if $w_{\mu}=w_{\nu}$ on the unit circle $S^1$. Real
 analytically, the Teichm\"uller space $T(X)$ is isomorphic to
 $\mathcal{B}^{-1,1}(\Del,\Gamma)/\sim$.

 Denote by $w^{\mu}$ the corresponding quasiconformal mapping if we
 extend
 $\mu\in\mathcal{B}^{-1,1}(\Del,\Gamma)$ to $\hat{\C}$ by setting it
 equal to zero outside $\Del$. $w^{\mu}$ is holomorphic on $\Del^*$
 and $\Gamma^{\mu}=w^{\mu}\circ\Gamma\circ (w^{\mu})^{-1}$ is no longer a
Fuchsian group, but a quasi--Fuchsian group. The corresponding
Riemann surface $X^{\mu}=\Gamma^{\mu}\bk w^{\mu}(\Del)$ is
biholomorphic to $X_{\mu}$. Complex analytically, the Teichmu\"uller
space is the quotient $\mathcal{B}^{-1,1}(\Del,\Gamma)/\sim$, where
$\mu\sim \nu$ if and only if $w^{\mu}=w^{\nu}$ on the unit circle.
For two compact Riemann surfaces $X$ and $Y$ having the same genus
$g$, their Teichm\"uller spaces $T(X)$ and $T(Y)$ are naturally
isomorphic and we use $T_g$ to denote the Teichm\"uller space of
compact Riemann surfaces of genus $g$.

 The tangent space and cotangent space at a
point $[\mu]$ of the Teichm\"uller space $T(X)$ can be naturally
identified with the space of harmonic Beltrami differentials and the
space of holomorphic quadratic differentials of $X^{\mu}$. The space
of holomorphic quadratic differentials of $X^{\mu}$ can be
identified with the space of $\Gamma_{\mu}$-automorphic $(2,0)$
differentials on $\Del$
--- $\Omega^{2,0}(\Del,\Gamma_{\mu})$, which consists of holomorphic
functions $q$ on $\Del$ satisfying
\begin{align*}
q(\gamma(z))\gamma'(z)^2= q(z),\hspace{1cm}\forall \;\gamma\in
\Gamma_{\mu}.
\end{align*}Correspondingly, the space of harmonic Beltrami
differentials of $X^{\mu}$ can be identified with
$\Omega^{-1,1}(\Del,\Gamma_{\mu})$, which is a subspace of
$\mathcal{A}^{-1,1}(\Del,\Gamma_{\mu})$ consisting of $\nu$ of the
form
$$\nu(z)=\rho(z)^{-1}\ov{q(z)}=\frac{(1-|z|^2)^2}{4}\ov{q(z)},$$where $q\in
\Omega^{2,0}(\Del,\Gamma_{\mu})$ and $\rho$ is the hyperbolic metric
density on $\Del$. The moduli space $\mathcal{M}(X)$ is the quotient
of the Teichm\"uller space $T(X)$ under the action of the mapping
class group. The Weil--Petersson metric on $T(X)$, which is defined
by
\begin{align*}
\langle \nu_{\alpha},\nu_{\beta}\rangle_{WP}
=\iint\limits_{\Gamma_{\mu}\bk\Del}\nu_{\alpha}(z)\overline{\nu_{\beta}(z)}\rho(z)
d^2z
\end{align*}on the tangent space of $T(X)$ at $[\mu]$, is modular
invariant and hence descends to a well-defined metric on the moduli
space $\mathcal{M}(X)$.

The universal Teichm\"uller space $T(1)$ can be defined similarly
with the group $\Gamma$ being the trivial group consisting of only
the identity element, i.e., $\Gamma=\{\text{id}\}$. More precisely,
$T(1)=\mathcal{B}^{-1,1}(\Del)/\sim$, where
$\mathcal{B}^{-1,1}(\Del)$ is the space of bounded functions on
$\Del$ with sup-norm less than one; and $\mu\sim\nu$ if and only if
$w_{\mu}\sim w_{\nu}$ on $S^1$. The cotangent space at any point of
$T(1)$ is naturally isomorphic to the Banach space
\begin{align*}A_{\infty}(\Del)=\left\{ q \;\text{holomorphic on
$\Del$}\;:\; \Vert q\Vert_{\infty}=\sup_{z\in \Del}
\rho(z)^{-1}|q(z)|<\infty\right\};
\end{align*}while the tangent space is identified with the Banach
space $$\Omega^{-1,1}(\Del)=\Bigl\{\rho^{-1}\bar{q}\;:\; q\in
A_{\infty}(\Del)\Bigr\}
$$ of harmonic Beltrami differentials on $\Del$. Obviously, the
 inner product
 \begin{align}\label{eq12_17_1}\langle \nu_{\alpha},\nu_{\beta}\rangle
  = \iint\limits_{ \Del}\nu_{\alpha}(z)\overline{\nu_{\beta}(z)}\rho(z)
d^2z\end{align}is not well--defined on $\Omega^{-1,1}(\Del)$. In
\cite{TT}, we showed that we can define a Hilbert structure on
$T(1)$ so that at any point, its tangent space is isomorphic to
\begin{align*}
H^{-1,1}(\Del)=\Bigl\{ \rho^{-1}\bar{q}\,:\, q\in A_2(\Del)\Bigr\},
\end{align*}where
\begin{align*}
A_2(\Del)=\left\{q \;\text{holomorphic on $\Del$}\;:\; \Vert
q\Vert_{2}^2=\iint\limits_{\Del}
|q(z)|^2\rho(z)^{-1}d^2z<\infty\right\}.
\end{align*}We denote the Teichm\"uller space with this Hilbert structure as $T_H(1)$.
The inner product \eqref{eq12_17_1} is well--defined on the tangent
space $H^{-1,1}(\Del)$ and we called the resulting metric on
$T_H(1)$ the Weil--Petersson metric.

Let $\Delta=-\rho^{-1}\pa\bar{\pa}$ be the Laplace--Beltrami
operator of the hyperbolic metric on $X$ and let
\begin{align}\label{eq12_19_1}G=\frac{1}{2}\left(\Delta+\frac{1}{2}\right)^{-1}\end{align} be one-half
of the resolvent of $\Delta$ at $\lambda=-1/2$. In \cite{Wolpert1},
Wolpert showed that the Riemann curvature tensor
$R_{\alpha\bar{\beta}\gamma\bar{\delta}}$ of the Weil--Petersson
metric at the tangent space of a point on the moduli space
correspond to the compact Riemann surface $X=\Gamma\bk\Del$ is given
by\footnote{Our convention differ with the convention of Wolpert by
a sign.}
\begin{align}\label{eq12_17_2}
R_{\alpha\bar{\beta}\lambda\bar{\delta}}=-\iint_{\Gamma\bk\Del}G(\nu_{\alpha}\bar{\nu}_{\beta})(
\nu_{\lambda}\bar{\nu}_{\delta})\rho
d^2z-\iint\limits_{\Gamma\bk\Del}G(\nu_\alpha\bar{\nu}_{\delta})(\nu_{\lambda}\bar{\nu}_{\beta})\rho
d^2z.
\end{align}In \cite{TT}, we generalized this result and showed that
formula \eqref{eq12_17_2} is still valid on $T_H(1)$ if
$\Gamma\bk\Del$ is replaced by $\Del$.

\section{Lower Bounds of curvatures of the Weil--Petersson metric on moduli space of compact Riemann surfaces}
In \cite{Wolpert1}, Wolpert has shown that the holomorphic
sectional, Ricci and scalar curvatures are bounded above by
$-\frac{1}{2\pi(g-1)}$, $-\frac{1}{2\pi(g-1)}$ and
$-\frac{3(3g-2)}{4\pi}$ respectively. These upper bounds depend on
the genus $g$. On the other hand, Huang showed that the sectional
curvature does not have a negative upper bound \cite{Huang2}.  Here
we would like to  find lower bounds for the curvatures, which only
depend on the injectivity radius of the corresponding Riemann
surface.

Recall that the holomorphic sectional and Ricci curvatures at a
point on the moduli space corresponding to the Riemann surface
$X=\Gamma\bk\Del$ in the direction  spanned by
$\nu_{\alpha}\in\Omega^{-1,1}(\Del,\Gamma)$ with
$\Vert\nu_{\alpha}\Vert_{WP}=1$, are given respectively by
\cite{Bochner, Wolpert1}:
\begin{align}\label{eq12_18_1}s_{\alpha}=R_{\alpha\bar{\alpha}\alpha\bar{\alpha}}=-2\iint\limits_{\Gamma\bk\Del}
G(|\nu_{\alpha}|^2)|\nu_{\alpha}|^2\rho d^2z,\end{align}and
\begin{align}\label{eq12_18_2}
\mathcal{R}_{\alpha\bar{\alpha}}=&\sum_{\beta=1}^{3g-3}
R_{\alpha\bar{\beta}\beta\bar{\alpha}}\\
=&-\sum_{\beta=1}^{3g-3}\left\{\iint\limits_{\Gamma\bk\Del}
G(\nu_{\alpha}\bar{\nu}_{\beta})\bar{\nu}_{\alpha}\nu_{\beta}\rho
d^2z +\iint\limits_{\Gamma\bk\Del}
G(|\nu_{\alpha}|^2)|\nu_{\beta}|^2\rho d^2z\right\},\nonumber
\end{align}where $\left\{\nu_1,\ldots,\nu_{3g-3}\right\}$ is an orthonormal basis of $\Omega^{-1,1}(\Del,\Gamma)$.
The scalar curvature $S$ is equal to the trace of the Ricci tensor:
\begin{align}\label{eq12_18_3}
S=\sum_{\alpha=1}^{3g-3}
\mathcal{R}_{\alpha\bar{\alpha}}=-\sum_{\alpha=1}^{3g-3}\sum_{\beta=1}^{3g-3}\left\{\iint\limits_{\Gamma\bk\Del}
G(\nu_{\alpha}\bar{\nu}_{\beta})\bar{\nu}_{\alpha}\nu_{\beta}\rho
d^2z +\iint\limits_{\Gamma\bk\Del}
G(|\nu_{\alpha}|^2)|\nu_{\beta}|^2\rho d^2z\right\}.
\end{align}
On the other hand, given two orthogonal tangent vectors
$\nu_{\alpha},\nu_{\beta}\in \Omega^{-1,1}(\Del,\Gamma)$ with
$\Vert\nu_{\alpha}\Vert_{WP}=\Vert\nu_{\beta}\Vert_{WP}=1$, the
sectional curvature of the plane spanned by the real tangent vectors
corresponding to $\nu_{\alpha}$ and $\nu_{\beta}$ is \cite{Bochner,
Wolpert1}
\begin{align}\label{eq12_18_4}
K_{\alpha,\beta}=&\frac{1}{4}\left(R_{\alpha\bar{\beta}\beta\bar{\alpha}}+R_{\beta\bar{\alpha}\alpha\bar{\beta}}
-R_{\alpha\bar{\beta}\alpha\bar{\beta}}-R_{\beta\bar{\alpha}\beta\bar{\alpha}}\right)\\
=&\text{Re}\iint\limits_{\Gamma\bk\Del}
G(\nu_{\alpha}\bar{\nu}_{\beta})\nu_{\alpha}\bar{\nu}_{\beta}\rho
d^2z -\frac{1}{2}\iint\limits_{\Gamma\bk\Del}
G(|\nu_{\alpha}|^2)|\nu_{\beta}|^2\rho d^2z
-\frac{1}{2}\iint\limits_{\Gamma\bk\Del}
G(\nu_{\alpha}\bar{\nu}_{\beta})\bar{\nu}_{\alpha}\nu_{\beta}\rho
d^2z.\nonumber
\end{align}
We have used the self--adjointness of $G$ to obtain the last
expression.

To obtain the lower bounds of curvatures, Huang \cite{Huang1} used
harmonic maps to  show that if a harmonic Beltrami differential has
unit Weil--Petersson norm, then its sup-norm is bounded above by a
constant depending on the injectivity radius of the underlying
Riemann surface. Here we reprove this result without resorting to
harmonic maps, which better reveals its elementary nature; and also
allows us to generalize this result to the universal Teichm\"uller
space later.

\begin{proposition}\label{proposition1}
Let $X=\Gamma\bk\Del$ be a compact Riemann surface with injectivity
radius $r_X$ and let $\nu\in\Omega^{-1,1}(\Del,\Gamma)$ be a
harmonic Beltrami differential of $X$. The ratio of the sup--norm of
$\nu$ to the Weil--Petersson norm of $\nu$ is bounded above by a
constant $C(r_X)$ depending only on $r_X$, i.e.,
\begin{align*}
\Vert\nu\Vert_{\infty} \leq C(r_X) \Vert \nu\Vert_{WP}.
\end{align*}The constant $C(r_X)$ can be chosen to be equal to
\begin{align}\label{eq12_17_4}
C(r_X)=\left\{
\frac{4\pi}{3}\left[1-\left(\frac{4e^{r_X}}{(e^{r_X}+1)^2}\right)^3\right]\right\}
^{-\frac{1}{2}}.
\end{align}
\end{proposition}
\begin{proof}
Let $z\in \Del$ and let $$\sigma_z(w) =\frac{z+w}{1+z\bar{w}},
\hspace{1cm}w\in\Del,$$ be a linear transformation preserving $\Del$
and mapping $0$ to $z$. Notice that $\nu\circ \sigma_z$ is a
harmonic Beltrami differential of the group
$\sigma^{-1}_z\circ\Gamma\circ\sigma_z$, and $\Vert \nu\circ
\sigma_{z}\Vert_{WP}=\Vert \nu\Vert_{WP}$, but $| \nu(z)|=|\nu\circ
\sigma_z(0)|$. Therefore, it suffices to verify that there exists a
constant $C(r_X)$ such that
$$|\nu(0)|\leq C(r_X) \Vert\nu\Vert_{WP}.$$

By definition, there exists $q\in \Omega^{2,0}(\Del,\Gamma)$ such
that $\nu = \rho^{-1} \bar{q}$. Being a holomorphic function on
$\Del$, $q$ has a Taylor series expansion on $\Del$ which can be
written as
\begin{align*}
q(z) =\sum_{n=2}^{\infty} (n^3-n) a_n z^{n-2}.
\end{align*}This implies that $$\nu(0)=\rho(0)^{-1}q(0)=\frac{3a_2}{2},$$ whereas
\begin{align*}
\Vert \nu\Vert_{WP}^2=\iint\limits_{\Gamma\bk\Del} |q(z)|^2
\rho(z)^{-1} d^2z.
\end{align*}By the definition of injectivity radius, we can choose a
fundamental domain $F$ for the action of $\Gamma$ on $\Del$ which
contains a hyperbolic disc $D(0, r)$ with center at $0$ and with
radius $r$, for any $r$ less than $r_X$. Elementary  hyperbolic
geometry gives us
\begin{align*}
D(0, r)=\left\{ z\in \C\,:\, |z|<\frac{e^r-1}{e^r+1}\right\}.
\end{align*}Therefore, for any $r\in (0, r_X)$, we have
\begin{align}\label{eq12_17_3}
\Vert \nu\Vert_{WP}^2=&\iint\limits_{F} |q(z)|^2 \rho(z)^{-1}
d^2z\nonumber\\\geq& \iint\limits_{D(0, r)}|q(z)|^2
\rho(z)^{-1}d^2z\nonumber\\
=&\int_0^{\frac{e^r-1}{e^r+1}}\int_0^{2\pi}\left|\sum_{n=2}^{\infty}(n^3-n)
a_n u^{n-2} e^{i(n-2)\theta}\right|^2\frac{(1-u^2)^2}{4} d\theta
udu\nonumber\\
=&2\pi\int_0^{\frac{e^r-1}{e^r+1}}\sum_{n=2}^{\infty} (n^3-n)^2
|a_n|^2
u^{2n-4} \frac{(1-u^2)^2}{4} udu\nonumber\\
\geq & 18\pi|a_2|^2\int_0^{\frac{e^r-1}{e^r+1}} (1-u^2)^2 udu\nonumber\\
=&\frac{4\pi}{3}|\nu(0)|^2
\left\{1-\left(\frac{4e^r}{(e^r+1)^2}\right)^3\right\}.
\end{align}Since this is true for all $r\in (0, r_X)$, we can
replace $r$ in \eqref{eq12_17_3} by $r_X$. Therefore, we have proved
the proposition with $C(r_X)$ equal to
\begin{align*}
C(r_X)=\left\{
\frac{4\pi}{3}\left[1-\left(\frac{4e^{r_X}}{(e^{r_X}+1)^2}\right)^3\right]\right\}
^{-\frac{1}{2}}.
\end{align*}
\end{proof}
Notice that $C(r_X)$ is a decreasing function of $r_X$. As $r_X$
approaches infinity, it approaches $\sqrt{3/4\pi}$. On the other
hand, as $r_X\rightarrow 0$, it behaves like
\begin{align}\label{eq12_19_5}
C(r_X)\sim \frac{1}{\sqrt{\pi} r_X}+O(1).
\end{align}

Before computing the lower bounds for the curvature, we  state a
useful lemma here.
\begin{lemma}\label{lemma1}Let $G$ be the positive self--adjoint
operator on $X$ defined by \eqref{eq12_19_1}.

\noindent A. For any  $f\in L^2(X,\R)$
$$\iint\limits_{\Gamma\bk\Del}
G(f)\bar{f}\rho d^2z\geq 0.$$

\noindent B. For any  $f\in L^2(X,\R)$,
$$\iint\limits_{\Gamma\bk\Del}
G(f)\rho d^2z =\iint\limits_{\Gamma\bk\Del} f\rho d^2z.$$

\noindent C. If  $f\in L^2(X,\R)$ is such that $f\geq 0$, then
$G(f)\geq 0$.

\noindent D. For any $f,g\in L^2(X,\R)$, $|G(fg)|\leq G(f^2)^{1/2}
G(g^2)^{1/2}$.
\end{lemma}
\begin{proof}
A is an immediate consequence of the positivity of $G$. B is proved
using the self-adjointness of $G$ and the fact that $G(1)=1$. C
follows from the fact that the kernel of $G$ --- $G(z,w)$, is a
positive function for all $z$ and $w$ (see \cite{TT}). D is the
Lemma 4.3 in \cite{Wolpert1}.
\end{proof}

\begin{proposition}\label{proposition2}
Let $X=\Gamma\bk\Del$ be a compact Riemann surface of genus $g$ with
injectivity radius $r_X$. At the point on the moduli space
$\mathcal{M}_g$ corresponding to $X$, the holomorphic sectional and
sectional curvatures of the Weil-Petersson metric is bounded below
by $-2C(r_X)^2$.

\end{proposition}
\begin{proof}The proof follows closely the proofs to obtain lower
and upper bounds given in \cite{Wolpert1, Huang1}. For completeness,
we repeat it  here. We consider the holomorphic sectional curvature
first. Given $\nu_{\alpha}\in\Omega^{-1,1}(\Del,\Gamma)$ with
$\Vert\nu_{\alpha}\Vert_{WP}=1$, we find from \eqref{eq12_18_1},
Proposition \ref{proposition1} and B and C of Lemma \ref{lemma1}
that
\begin{align*}
-s_{\alpha}\leq &
2C(r_X)^2\iint\limits_{\Gamma\bk\Del}G(|\nu_{\alpha}|^2)\rho d^2z\\
=&2C(r_X)^2\iint\limits_{\Gamma\bk\Del} |\nu_{\alpha}|^2\rho
d^2z=2C(r_X)^2.
\end{align*}This proves the statement for holomorphic sectional
curvature. For the sectional curvature, we use formula
\eqref{eq12_18_4}. Notice that Cauchy--Schwarz inequality, the
positivity of the kernel of $G$, and D of Lemma \ref{lemma1} give us
\begin{align*}
\left|\iint\limits_{\Gamma\bk\Del}
G(\nu_{\alpha}\bar{\nu}_{\beta})\nu_{\alpha}\bar{\nu}_{\beta}\rho
d^2z\right|\leq & \iint\limits_{\Gamma\bk\Del} G(|\nu_{\alpha}||
\nu_{\beta}|)|\nu_{\alpha}||\nu_{\beta}|\rho d^2z\\
\leq &\iint\limits_{\Gamma\bk\Del}  G(|\nu_{\alpha}|^2)^{1/2}
G(|\nu_{\beta}|^2)^{1/2} |\nu_{\alpha}||\nu_{\beta}|\rho d^2z\\ \leq
&
\left\{\iint\limits_{\Gamma\bk\Del}G(|\nu_{\alpha}|^2)|\nu_{\beta}|^2
\rho
d^2z\right\}^{1/2}\left\{\iint\limits_{\Gamma\bk\Del}G(|\nu_{\beta}|^2)|\nu_{\alpha}|^2
\rho d^2z\right\}^{1/2}\\
=&\iint\limits_{\Gamma\bk\Del}G(|\nu_{\alpha}|^2)|\nu_{\beta}|^2
\rho d^2z.
\end{align*}Similarly,
\begin{align}\label{eq12_19_3}
\left|\iint\limits_{\Gamma\bk\Del}
G(\nu_{\alpha}\bar{\nu}_{\beta})\bar{\nu}_{\alpha}\nu_{\beta}\rho
d^2z\right|\leq\iint\limits_{\Gamma\bk\Del}G(|\nu_{\alpha}|^2)|\nu_{\beta}|^2
\rho d^2z.
\end{align}Therefore,
\begin{align*}
K_{\alpha,\beta}\geq -2\iint\limits_{\Gamma\bk\Del}
G(|\nu_{\alpha}|^2)|\nu_{\beta}|^2\rho d^2z.
\end{align*}The same reasoning as in the case of holomorphic
sectional curvature shows that $$K_{\alpha,\beta}\geq -2C(r_X)^2.$$

\end{proof}
If we naively use the approach above to find the lower bounds for
the Ricci and scalar curvatures, we will find that the Ricci and
scalar curvatures are bounded below by $-2(3g-3)C(r_X)^2$ and
$-2(3g-3)^2C(r_X)^2$ respectively. In fact,  these bounds are
obtained in \cite{Huang1}. However, by doing slightly more work, we
can greatly improve the bounds. Observe that given an orthonormal
basis $\{\nu_1,\ldots, \nu_{3g-3}\}$ of
$\Omega^{-1,1}(\Del,\Gamma)$, the
kernel\begin{align}\label{eq4_11_1}
P(z,w)=\rho(z)\rho(w)\sum_{\beta=1}^{3g-3}
\ov{\nu_{\beta}(z)}\nu_{\beta}(w)
\end{align}is the kernel of the projection operator mapping bounded quadratic
differentials to holomorphic quadratic differentials. Namely, for
any bounded quadratic differential $q$ of $X$,
\begin{align*}
(Pq)(z):= \iint\limits_{\Gamma\bk\Del} P(z,w) q(w) \rho(w)^{-1}d^2w
\end{align*}is a holomorphic quadratic differential, and $Pq=q$ if
and only if $q$ is holomorphic.  Let
\begin{align}\label{eq4_11_2}\Lambda = \sup_{z\in
\Del}\sum_{\beta=1}^{3g-3}|\nu_{\beta}(z)|^2=\sup_{z\in\Del}
\rho(z)^{-2}P(z,z).\end{align}Obviously,
\begin{align}\label{eq4_11_3}
\Lambda \leq \sup_{z\in\Del}\sup_{w\in \Del}
\left|\sum_{\beta=1}^{3g-3}\ov{\nu_{\beta}(z)}\nu_{\beta}(w)\right|.
\end{align}For fixed $z$,
$\sum_{\beta=1}^{3g-3}\ov{\nu_{\beta}(z)}\nu_{\beta}(w)\in
\Omega^{-1,1}(\Del,\Gamma)$. By \eqref{eq4_11_1}, its
Weil--Petersson norm is
\begin{align*}
\left\Vert
\sum_{\beta=1}^{3g-3}\ov{\nu_{\beta}(z)}\nu_{\beta}(w)\right\Vert_{WP}^2=&\rho(z)^{-2}
\iint\limits_{\Gamma\bk\Del}
P(z,w)\overline{P(z,w)}\rho(w)^{-1}d^2w.\end{align*}Since as a
function of $w$, $\overline{P(z,w)}=P(w,z)$ is a holomorphic
quadratic differential, the projection property of $P$ implies that
\begin{align*}
\iint\limits_{\Gamma\bk\Del}
P(z,w)\overline{P(z,w)}\rho(w)^{-1}d^2w=&\overline{P(z,z)}=P(z,z).
\end{align*}Therefore,
\begin{align*}
\left\Vert
\sum_{\beta=1}^{3g-3}\ov{\nu_{\beta}(z)}\nu_{\beta}(w)\right\Vert_{WP}^2=\rho(z)^{-2}
P(z,z).
\end{align*}Using
Proposition \ref{proposition1}, this gives
\begin{align*}
\sup_{w\in \Del}
\left|\sum_{\beta=1}^{3g-3}\ov{\nu_{\beta}(z)}\nu_{\beta}(w)\right|\leq
C(r_X) \sqrt{\rho(z)^{-2}P(z,z)}.
\end{align*}Consequently, \eqref{eq4_11_3} and \eqref{eq4_11_2}
imply that
\begin{align*}
\Lambda \leq \sup_{z\in\Del}\sup_{w\in \Del}
\left|\sum_{\beta=1}^{3g-3}\ov{\nu_{\beta}(z)}\nu_{\beta}(w)\right|\leq
C(r_X) \sup_{z\in \Del}
\sqrt{\rho(z)^{-2}P(z,z)}=C(r_X)\Lambda^{1/2}.
\end{align*}
In other words,\begin{align}\label{eq12_19_4}\Lambda=\sup_{z\in
\Del}\sum_{\beta=1}^{3g-3}|\nu_{\beta}(z)|^2\leq
C(r_X)^2.\end{align} Notice that we greatly improve the naive bound
$\Lambda\leq (3g-3)C(r_X)^2$ to $\Lambda\leq C(r_X)^2$. Now we can
prove the following lower bounds for Ricci and scalar curvatures:
\begin{proposition}\label{proposition3}~

\begin{itemize}
\item[A.] The Ricci curvature of the Weil--Petersson metric is bounded
below by $-2C(r_X)^2$.\\
\item[B.]  The scalar curvature of the Weil--Petersson metric is
bounded below by $-2(3g-3) C(r_X)^2$.\end{itemize}
\end{proposition}
\begin{proof}  Using  \eqref{eq12_18_2}, \eqref{eq12_18_3}
and \eqref{eq12_19_3}, we have
\begin{align*}
\mathcal{R}_{\alpha\bar{\alpha}}\geq &
-2\sum_{\beta=1}^{3g-3}\int\limits_{\Gamma\bk\Del}
G(|\nu_{\alpha}|^2)|\nu_{\beta}|^2\rho d^2z,\\
S\geq &
-2\sum_{\alpha=1}^{3g-3}\sum_{\beta=1}^{3g-3}\int\limits_{\Gamma\bk\Del}
G(|\nu_{\alpha}|^2)|\nu_{\beta}|^2\rho d^2z.
\end{align*}Eq. \eqref{eq12_19_4} and the same method used in the proof
of Proposition \ref{proposition2} then give us immediately
\begin{align*}
\mathcal{R}_{\alpha\bar{\alpha}}\geq -2C(r_X)^2, \hspace{1cm} S\geq
-2(3g-3) C(r_X)^2.
\end{align*}\end{proof}
Notice that we have established that the Ricci curvature of the
Weil--Petersson metric is bounded below by a constant depending only
on the injectivity radius of the corresponding Riemann surface, but
independent of the genus. This substantially improves the result of
\cite{Huang1}.

We would also like to remark that tending to the boundary of the
moduli spaces, the injectivity radius $r_X=inj(X)$ decreases to
zero. The results of Propositions \ref{proposition2} and
\ref{proposition3} and the estimate \eqref{eq12_19_5} show that when
$r_X\rightarrow 0$, the holomorphic sectional, sectional and Ricci
curvatures of the Weil--Petersson metric are all bounded below by a
constant of order $1/r_X^2$. It is interesting to compare this with
the asymptotics of the curvatures obtained in \cite{Schumacher}.

\section{Bounds of curvatures on the universal Teichm\"uller space}

In \cite{TT}, we have shown that the holomorphic sectional and
sectional curvatures of the Weil--Petersson metric on the Hilbert
manifold $T_H(1)$ are negative. We also showed that $T_H(1)$ is a
Kahler--Einstein manifold with constant Ricci curvature
$-\frac{13}{12\pi}$. In this section, we show that the holomorphic
sectional and sectional curvatures are  bounded below by a universal
constant. We also show that these curvatures do not have negative
upper bounds.

An analog of Proposition \ref{proposition1} for the universal
Teichm\"uller space $T_H(1)$ is
\begin{lemma}\label{lemma2}
Let $\nu\in \Omega^{-1,1}(\Del)$ be a harmonic Beltrami differential
on $\Del$. Then
\begin{align}\label{eq12_19_6}
\Vert \nu\Vert_{\infty}\leq \sqrt{\frac{3}{4\pi}}\Vert\nu\Vert_{WP}.
\end{align}
\end{lemma}This can be considered as the limiting case of
Proposition \ref{proposition1} when $r_X\rightarrow \infty$. In
fact, in the present situation, the corresponding Riemann surface is
isomorphic to the disc which has infinite hyperbolic radius. Another
proof of \eqref{eq12_19_6} is given in the proof of Lemma 2.1 in
\cite{TT}.

Using Lemma \ref{lemma2}, one obtains immediately as in Proposition
\ref{proposition2} that
\begin{proposition}
On the universal Teichm\"uller space $T_H(1)$, the holomorphic
sectional and sectional curvatures are bounded below by
$-\frac{3}{2\pi}$.
\end{proposition}

To prove the statements about upper bounds, we define for $n\geq 2$,
$$\nu_n=\rho(z)^{-1}\sqrt{\frac{2(n^3-n)}{\pi}}\bar{z}^{n-2}.$$
It is easy to show that $\{\nu_2,\nu_3,\ldots\}$ is an orthonormal
basis of $H^{-1,1}(\Del)$. It is elementary to find the sup--norm of
$\nu_n$ explicitly:
\begin{lemma}
For $n\geq 2$, $$\Vert
\nu_n\Vert_{\infty}=\sqrt{\frac{2(n^3-n)}{\pi}}\frac{4}{(n+2)^2}\left(\frac{n-2}{n+2}\right)^{\frac{n-2}{2}}.$$
\end{lemma}\begin{proof}
For $n\geq 2$, define
\begin{align*}
h_n(r) = (1-r^2)^2 r^{n-2},\hspace{1cm}r\in[0,1].
\end{align*}Then $h_n$ is a nonnegative function and
\begin{align*}
h_n'(r)=r^{n-3}(1-r^2)((n-2)-(n+2)r^2).
\end{align*}This implies that $h_n(r)$ has maximum at
$r=\sqrt{(n-2)/(n+2)}$ and its maximum value is
\begin{align*}
\max_{r\in[0,1]}
h_n(r)=\frac{16}{(n+2)^2}\left(\frac{n-2}{n+2}\right)^{\frac{n-2}{2}}.
\end{align*}The assertion follows.
\end{proof}
Notice that $\Vert\nu_2\Vert_{\infty}=\sqrt{3/(4\pi)}$. This shows
that the result of Lemma \ref{lemma2} is sharp. On the other hand,
it is easy to see that
\begin{align}\label{eq12_20_1}
\Vert \nu_n\Vert_{\infty}\leq \sqrt{\frac{32}{\pi(n+2)}}\rightarrow
0\hspace{1cm}\text{as}\;\; n\rightarrow \infty.
\end{align}Using this, we can prove that
\begin{proposition} On the universal Teichm\"uller space
$T_H(1)$, the holomorphic sectional and sectional curvatures do not
have negative upper bounds.
\end{proposition} \begin{proof}
For the holomorphic sectional curvature, we obtain as in the proof
of Proposition \ref{proposition2} that
\begin{align*}
|s_n|=2\iint\limits_{\Del}G(|\nu_n|^2)|\nu_n|^2\rho d^2z\leq
2\Vert\nu_n\Vert_{\infty}^2.
\end{align*}
On the other hand, the proof of Proposition \ref{proposition2} shows
that the sectional curvature $K_{m,n}$ \eqref{eq12_18_4} is bounded
by
\begin{align*}
|K_{m,n}|\leq 2\iint\limits_{\Del}G(|\nu_m|^2) |\nu_n|^2\rho
d^2z\leq 2\Vert\nu_n\Vert_{\infty}^2.
\end{align*}Since by \eqref{eq12_20_1},
 $\Vert\nu_n\Vert_{\infty}\rightarrow 0$ as $n\rightarrow \infty$, we conclude that the
holomorphic sectional and sectional curvatures do not have negative
upper bounds.
\end{proof}

 \vspace{1cm} \noindent \textbf{Acknowledgement}\;
The author would like to thank the Ministry of Science, Technology
and Innovation of Malaysia for funding this project under
eScienceFund 06-02-01-SF0021.


\begin{thebibliography}{10}


\bibitem{Ahlfors}
Lars~V. Ahlfors, \emph{Curvature properties of {T}eichm\"uller's
space}, J.
  Analyse Math. \textbf{9} (1961/1962), 161--176.

\bibitem{Bochner}
S.~Bochner, \emph{Curvature in {H}ermitian metric}, Bull. Amer.
Math. Soc.
  \textbf{53} (1947), 179--195.

\bibitem{Huang2}
Zheng Huang, \emph{Asymptotic flatness of the {W}eil-{P}etersson
metric on
  {T}eichm\"uller space}, Geom. Dedicata \textbf{110} (2005), 81--102.

\bibitem{Huang1}
Zheng Huang, \emph{The {W}eil-{P}etersson geometry on the thick part
of the moduli
  space of {R}iemann surfaces}, Proc. Amer. Math. Soc. \textbf{135} (2007),
  no.~10, 3309--3316 (electronic).

\bibitem{Huang3}
Zheng Huang, \emph{On asymptotic Weil--Petersson geometry of
Teichm\"uller space of Riemann surface}, preprint arXiv:
math.DG/0405228, to appear in Asian J. Math.

\bibitem{Jost}
J{\"u}rgen Jost, \emph{Harmonic maps and curvature computations in
  {T}eichm\"uller theory}, Ann. Acad. Sci. Fenn. Ser. A I Math. \textbf{16}
  (1991), no.~1, 13--46.



\bibitem{Liu_Sun_Yau_1}
Kefeng Liu, Xiaofeng Sun, and Shing-Tung Yau, \emph{Canonical
metrics on the
  moduli space of {R}iemann surfaces. {I}}, J. Differential Geom. \textbf{68}
  (2004), no.~3, 571--637.

\bibitem{Liu_Sun_Yau_2}
Kefeng Liu, Xiaofeng Sun, and Shing-Tung Yau, \emph{Canonical
metrics on the moduli space of {R}iemann surfaces.
  {II}}, J. Differential Geom. \textbf{69} (2005), no.~1, 163--216.

\bibitem{McMullen}
Curtis~T. McMullen, \emph{The moduli space of {R}iemann surfaces is
{K}\"ahler
  hyperbolic}, Ann. of Math. (2) \textbf{151} (2000), no.~1, 327--357.

\bibitem{Royden}
H.~L. Royden, \emph{Intrinsic metrics on {T}eichm\"uller space},
Proceedings of
  the International Congress of Mathematicians (Vancouver, B. C., 1974), Vol.
  2, Canad. Math. Congress, Montreal, Que., 1975, pp.~217--221.

\bibitem{Schumacher}
Georg Schumacher, \emph{Harmonic maps of the moduli space of compact
{R}iemann
  surfaces}, Math. Ann. \textbf{275} (1986), no.~3, 455--466.

\bibitem{TT}
Leon~A. Takhtajan and Lee-Peng Teo, \emph{Weil-{P}etersson metric on
the
  universal {T}eichm\"uller space}, Mem. Amer. Math. Soc. \textbf{183} (2006),
  no.~861, viii+119.

\bibitem{Trapani}
Stefano Trapani, \emph{On the determinant of the bundle of
meromorphic
  quadratic differentials on the {D}eligne-{M}umford compactification of the
  moduli space of {R}iemann surfaces}, Math. Ann. \textbf{293} (1992), no.~4,
  681--705.

\bibitem{Tromba}
A.~J. Tromba, \emph{On a natural algebraic affine connection on the
space of
  almost complex structures and the curvature of {T}eichm\"uller space with
  respect to its {W}eil-{P}etersson metric}, Manuscripta Math. \textbf{56}
  (1986), no.~4, 475--497.

\bibitem{Wilson}
P.~M.~H. Wilson, \emph{Sectional curvatures of {K}\"ahler moduli},
Math. Ann.
  \textbf{330} (2004), no.~4, 631--664.

\bibitem{Wolpert1}
Scott~A. Wolpert, \emph{Chern forms and the {R}iemann tensor for the
moduli
  space of curves}, Invent. Math. \textbf{85} (1986), no.~1, 119--145.

\bibitem{Wolpert2}
Scott~A. Wolpert, \emph{Geometry of the {W}eil-{P}etersson
completion of {T}eichm\"uller
  space}, Surveys in differential geometry, Vol.\ VIII (Boston, MA, 2002),
  Surv. Differ. Geom., VIII, Int. Press, Somerville, MA, 2003, pp.~357--393.

\bibitem{Yamada}
Sumio Yamada, \emph{On the geometry of {W}eil-{P}etersson completion
of
  {T}eichm\"uller spaces}, Math. Res. Lett. \textbf{11} (2004), no.~2-3,
  327--344.




\end{thebibliography}
\end{document}